\documentclass{article}
\usepackage{amssymb}
\usepackage{latexsym}
\usepackage{amsmath}
\usepackage{mathrsfs}
\usepackage[all]{xy}
\usepackage{enumerate}
\usepackage{amsthm}
\usepackage[dvips]{graphicx}
\usepackage{psfrag}
\usepackage{ifthen}

\newtheorem{theorem}{Theorem}[section]
\newtheorem{lemma}[theorem]{Lemma}
\newtheorem{metalemma}[theorem]{Meta Lemma}
\newtheorem{proposition}[theorem]{Proposition}
\newtheorem{corollary}[theorem]{Corollary}
\newtheorem{problem}[theorem]{Problem}

\theoremstyle{definition}
\newtheorem{definition}[theorem]{Definition}

\newtheorem{example}[theorem]{Example}

\newenvironment{thm}{\begin{theorem}}{
\end{theorem}}

\newenvironment{prop}{\begin{proposition}}{
\end{proposition}}

\newenvironment{cor}{\begin{corollary}}{
\end{corollary}}

\newenvironment{lem}{\begin{lemma}}{
\end{lemma}}

\newenvironment{defn}{\begin{definition}}{
\end{definition}}

\newcommand{\comment}[1]
	   {\ifthenelse{\equal{\showcomments}{yes}}
	     {\footnotemark\marginpar{\sffamily{\tiny
		   \addtocounter{footnote}{-1}\footnotemark#1
}\normalfont}}{}}

\newcommand{\noi}{\noindent}
\newcommand{\tr}{\textrm}
\newcommand{\showcomments}{yes}

\newcommand{\np}{{NP}}
\newcommand{\pmo}{{\pm 1}}
\newcommand{\mpo}{{\mp 1}}
\newcommand{\mo}{{-1}}
\newcommand{\bel}[1]{\begin{equation}\label{#1}}
\newcommand{\ee}{\end{equation}}

\newcommand{\bdy}{\partial}

\newcommand{\ol}{\overline}
\newcommand{\E}{\mathbb{E}}

\newcommand{\LBA}{\left\{\begin{array}}
\newcommand{\EAR}{\end{array}\right.}
\title{Quadratic equations over free groups are NP-complete}
\author{O. Kharlampovich, I.G. Lys\"enok, A.G. Myasnikov, N.W.M.
  Touikan}
\begin{document}
\maketitle
\abstract{We prove that the problems of deciding whether a quadratic
  equation over a free group has a solution is \np-complete}
\section{Introduction}
The study of quadratic equations over free groups probably started
with \cite{Malcev-1962} and has been deepened extensively ever since.
One of the reasons research in this topic has been so fruitful is the
deep connection between quadratic equations and the topology of
surfaces. 

In this paper we will show that the problem of deciding if a quadratic
equation over a free group is NP-complete. This problem was shown to be
decidable in \cite{Comerford-Edmunds-1981}. In addition it was shown in
\cite{Olshanskii-1989}, \cite{Grigorchuk-Kurchanov-1989}, and
\cite{Grigorchuk-Lysenok-1992} that if $n$, the number of variables,
is fixed, then deciding if a quadratic equation has a solution can be
done in time polynomial in the sum of the lengths of the coefficients.
These results imply that the problem is solvable in at most
exponential time. We will improve on this by proving (see Theorem
\ref{thm:in-np}) that deciding if an quadratic equation over a free
group has a solution is in NP.

In \cite{Diekert-Robson-1999} it is shown that deciding if a quadratic
word equation has a solution is NP-hard. We will prove (see Theorem
\ref{thm:np-hard}) that deciding if a quadratic equations over a free
group has a solution is also NP-hard. Our proofs are geometric,
relying on the topological results of \cite{Olshanskii-1989} and disc
diagram techniques.

\section{Quadratic equations over free groups are in NP}\label{sec:in-np}
For a finite alphabet alphabet $S$ we denote by $S^*$ the free monoid
with involution with basis $S$ and for $w \in S^*$, we denote by
$w^\mo$ its involution. We denote by $F(S)$ the free group on $S$.
\subsection{Standard form}
A quadratic equation $E$ with variables $\{x_i,y_i,z_j\}$ and
coefficients $\{w_i,d\} \in F(A)$ is said to be in \emph{standard
  form} if its coefficients are expressed as freely and cyclically
reduced words in $A^*$ and $E$ has either the form:
\bel{eqn:orientable} \prod_{i=1}^g [x_i,y_i] \prod_{j=1}^{m-1} z_j^\mo
w_j z_j d = 1 \tr{~or~} \prod_{i=1}^g [x_i,y_i]=1\ee where
$[x,y]=x^\mo y^\mo xy$, in which case we say it is \emph{orientable}
or it has the form \bel{eqn:non-orientable} \prod_{i=1}^g x_i^2
\prod_{j=1}^{m-1} z_j^\mo w_j z_j d = 1 \tr{~or~} \prod_{i=1}^g
x_i^2d=1\ee in which case we say it is non orientable.  The
\emph{genus} if a quadratic equation is the number $g$ in
(\ref{eqn:orientable}) and (\ref{eqn:non-orientable}) and $m$ is the
number of coefficients. If $g=0$ then we will define $E$ to be
orientable. If $E$ is a quadratic equation we define its \emph{reduced
  euler characteristic}, $\ol{\chi}$ as follows:
\[\ol{\chi}(E) =\left\{
  \begin{array}{l}2-2g \tr{~if~} E \tr{~is orientable} \\ 2-g
    \tr{~if~} E \tr{~is not orientable}\end{array}\right.\] It is a
well known fact that an arbitrary quadratic equation over a free group
can be brought to a standard form in time polynomial in its length.

\subsection{Ol'shanskii's result}\label{sec:olshanskii}
In sections 2.3 and 2.4 of \cite{Olshanskii-1989} it is shown that a
quadratic equation $E$ in standard form has a solution if and only if
for some $n \leq 3(m-\ol{\chi}(E))$,
\begin{itemize}
\item[(i)] there is a set $P = \{p_1,\ldots p_n\}$ of variables and a
  collection of $m$ discs $D_1,\ldots, D_m$ such that,
\item[(ii)] the boundaries of these discs are directed labeled graphs
  such that each edge has a label in $P$ and each $p_j \in P$ occurs
  exactly twice in the union of boundaries;
\item[(iii)] if we glue the discs together by edges with the same
  label, respecting the edge orientations, then we will have a
  collection $\Sigma_0,\ldots,\Sigma_l$ of closed surfaces and the
  following inequalites: if $E$ is orientable then each $\Sigma_i$ is
  orientable and \[\biggl(\sum_{i=0}^{l} \chi(\Sigma_i)\biggr) - 2l
  \geq \ol{\chi}(E)\] if $E$ is non-orientable either at least one
  $\Sigma_i$ is non-orientable and \[\biggl(\sum_{i=0}^{l}
  \chi(\Sigma_i)\biggr) - 2l \geq \ol{\chi}(E)\] or, each $\Sigma_i$
  is orientable and
  \[\biggl(\sum_{i=0}^{l} \chi(\Sigma_i)\biggr) - 2l \geq \ol{\chi}(E) +2 \] and
\item[(iv)] there is a mapping $P \rightarrow A^*$ such that upon
  substitution, the coefficients $w_1,\ldots,w_{m-1}$ and $d$ can be
  read without cancellations around the boundaries of
  $D_1,\ldots,D_{m-1}$ and $D_m$, respectively; and finally that
\item[(v)] if $E$ is orientable the discs $D_1,\ldots,D_m$ can be
  oriented so that $w_i$ is read clockwise around $\bdy D_i$ and $d$
  is read clockwise around $\bdy D_m$, moreover all these orientations
  must be compatible with the glueings.
\end{itemize} We note that the bounds in (iii) are not given
explicitly in that paper, but they follow immediately from the
discussion on cutting up diagrams into so-called simple diagrams, see
\cite{Olshanskii-1989} for details.

\subsection{The certificate}
The result of section \ref{sec:olshanskii} enables us to construct a
good certificate.
\begin{thm}\label{thm:in-np} For a quadratic equation $E$ in standard
  form there is a certificate of size bounded by
  $2(|w_1|+\ldots+|w_m|+|d|+3(2g+m))$ that can be checked in
  polynomial time.
\end{thm}
\begin{proof} The certificate will consist of the
  following:\begin{enumerate}
  \item A collection of variables $P=\{p_1,\ldots,p_n\}$.
  \item A collection of substitutions $\ol{\psi}=\{p_i \mapsto a_i,
    i=1\ldots n\}$ where $a_i \in A^*$ and $n < 3(2g+m)$.  
  \item A collection of words in $P^*$
    \[\mathcal{C}=\LBA{l}C_1=p_{11}^{\epsilon_{11}} \ldots
    p_{11}^{\epsilon_{1j_1}}\\ \ldots\\ C_m= p_{m1}^{\epsilon_{m1}}
    \ldots p_{mj_m}^{\epsilon_{mj_m}}\EAR\] with $p_{ij} \in P,
    \epsilon_{ij}\in \{-1,1\}$ and each $p_i \in P$ occuring exactly
    twice.
\end{enumerate}
The $C_i's$ are represent the labels of the boundaries of the discs
$D_1,\ldots D_l$ so checking contitions (i) and (ii) of Section
\ref{sec:olshanskii} can be done quickly, moreover we see that the
size of $\mathcal{C}$ is at most $2n \leq 6(2g+m)$. 

$\ol{\psi}$ extends to a monoid homomorphism $\psi: P^* \rightarrow
A^*$. (iv) can also be verified quickly since for $i=1,\ldots m-1$ we
just need to check that some cyclic permutation of $\psi(C_i)$ is
equal to $w_i$ and some cyclic permutation of $\psi(C_m)$ is equal do
$d$. Moreover, since the equality is graphical we have that
\[|a_1| + \ldots |a_n| \leq |w_1|+\ldots+|w_m|+|d|\] Therefore the
size of the certificate is bounded as advertised. All that is left is
to determine the topology of the glued together discs.  We describe
the algorithm without too much detail.
\\ \\
\noi {\bf Step 1:} \emph{Built a forest of discs:} We make a graph
$\Gamma$ such that each vertex $v_i \in V(\Gamma)$ corresonds the disc
$D_i$ and each edge $e_j\in E(\Gamma)$ corresponds to the variable $p_j
\in P$. The edge $e_k$ goes from $v_i$ to $v_j$ if and only if the
variable $p_k$ occurs in the boudary of $D_i$ and in the boundary of
$D_j$ or if $i=j$ then there are two different occurences of the
variable $p_k$. We construct a spanning forest $\mathcal{F}$. This
enables us to count the number of connected components
$\Sigma_0,\ldots,\Sigma_l$.

\noi {\bf Step 2:} \emph{Determine orientability:} For each maximal
tree $T_r \subset \mathcal{F}$ we get a ``tree of discs'' by glueing
together only the pairs of edges whose labels correspond to elements
of $E(T_r)$. The resulting tree of discs is a simply connected
topological space that can be embedded in the plane and we can read a
cyclic word $c(T_r)$ in $P^*$ along its boundary. The surface $\Sigma_r$
obtained by glueing together the remaining paired edges of the tree of
discs will be orientable only if whenever $p_j^\pmo$ occurs in $c(T_r)$
then $p_j^\mpo$ also occurs. We can also check (v) at this point.

\noi {\bf Step 3:} \emph{Compute Euler characteristic:} The
identification of the boundary of the discs with graphs, enables us to
think of the discs as polygons. If a disc $D_i$ has $N_i$ sides then we
give each corner of $D_i$ an angle of $\pi(N_i - 2)/N_i$. Then for each
tree of discs produced in the previous step, we identify the remaining
pairs of edges to get the surfaces $\Sigma_0,\ldots \Sigma_l$, which
now have an extra angular structure. To each $\Sigma_i$, we can apply
the Combinatorial Gauss-Bonnet Theorem which states that for an angled
two-complex $X$,\[ 2\pi \chi(X) = \sum_{f\in X^{(2)}} \kappa(f) +
\sum_{v \in X^{(0)}} \kappa(v)\] where $X^{(2)}$ is the set of faces
and $X^{(0)}$ is the set of vertices. This angle assingment gives each
face $f$ a curvature $\kappa(f)=0$ and each vertex has curvature
\[\kappa(v)=2\pi - \bigl(\sum_{c \in link(v)}
\measuredangle(c)\bigl)\] i.e. $\kappa(v)$ is $2\pi$ minus the sum of
the angles that meet at $v$.

With an appropriate data structure one can perform steps 1-3 (not
necessarily in sequential order) in at most quadratic time in the size
of $\mathcal{C}$. Once all that is done, verifying the inequalities of
$(iii)$ is easy and we are finished.
\end{proof}

\section{Quadratic equations over free groups are NP-hard}\label{sec:np-hard}
We will present the bin packing problem which is known to be
NP-complete and show that it is equivalent deciding if a certain type
of quadratic equation has a solution.

\subsection{Bin Packing}
\begin{problem}[Bin Packing]\label{prob:bin-packing}~
\begin{itemize}
\item INPUT: A $k-$tuple of positive integers $(r_1,\ldots,r_k)$ and
  positive integers $B, N$.
\item QUESTION: Is there a partition of $\{1,\ldots,k\}$ into $N$
  subsets \[ \{1,\ldots,k\} = S_1 \sqcup \ldots \sqcup S_N\] such that
  for each $i=1,\ldots, N$ we have \bel{eqn:bin}\sum_{j \in S_i} r_j
  \leq B\ee
\end{itemize}
\end{problem}

This problem is NP-hard in the strong sense (see \cite{Garey-Johnson}
p.226), i.e. there are NP-hard instances of this problem when both $B$
and the $r_j$ are bounded by a polynomial function of of $k$.

Let $t=NB-\sum_{i=1}^{k}r_i$. Then by replacing $(r_1,\ldots,r_k)$ by
the $k+t$-tuple $(r_1,\ldots,r_k,\ldots,1,\ldots,1)$ we can assume
that the inequalities (\ref{eqn:bin}) are actually equalities. This
modified version is still $NP$ hard in the strong sense. We state it
explicitly:
\begin{problem}[Exact Bin Packing]\label{prob:mod-bin-packing}~
\begin{itemize}
\item INPUT: A $k-$tuple of positive integers $(r_1,\ldots,r_k)$ and
  positive integers $B, N$.
\item QUESTION: Is there a partition of $\{1,\ldots,k\}$ into $N$
  subsets \[ \{1,\ldots,k\} = S_1 \sqcup \ldots \sqcup S_N\] such that
  for each $i=1,\ldots, N$ we have \bel{eqn:mod-bin}\sum_{j \in S_i} r_j =
  B\ee
\end{itemize}
\end{problem}
The authors warmly thank Laszlo Babai for drawing their attention to
this problem in connection to tiling problems.
\subsection{Tiling discs} Throughout this section we will consider the
discs to be embedded in $\E^2$ and will always read clockwise around
closed curves.
\begin{defn} A \emph{$[a,b^n]$-disc} is a disc as in section
  \ref{sec:olshanskii} along whose boundary one can read the cyclic
  word $[a,b^n]$.
\end{defn}

\begin{defn} A \emph{$[a,b^n]$-ribbon} is a rectangular cell complex obtained by
  attaching $[a,b^j]$-discs by their $a$-labeled edges, such that we
  can read $[a,b^n]$ along its boundary. The \emph{top} of an
  $[a,b^n]$ ribbon is the boundary subpath along which we can read the
  word $b^{-n}$, the \emph{bottom} is the boundary subpath along which
  we can read the word $b^{n}$.
\end{defn}

\begin{defn} Let $D$ be a disc tiled by $[a,b^n]$-discs, we define the
  \emph{$a-$pattern} of $D$ to be a graph defined as
  follows:\begin{enumerate}
  \item In the middle of each $a$-labeled edge put a vertex.
  \item Between any two vertices contained in the same $[a,b^n]$-disc
    draw an edge.
  \end{enumerate}
  Connected components of $a-$patterns are called \emph{$a-$tracks}
\end{defn}

\begin{lem}\label{lem:circ} A disc $D$ tiled by finitely many
  $[a,b^n]$-discs cannot have any circular $a-$tracks.
\end{lem}
\begin{proof}
It is clear that every $a-$track is a graph whose vertices have
valency at most 2. If an $a-$track $t$ has vertices of valency 1 then
they must lie on $\bdy D$. 

Suppose towards a contradiction that $D$ has a circular $a$-track
$c$. Then $c$ divides $D$ into two components: an interior and an
exterior. If we examine the interior we see that it is a planar union
of discs with only the letter $b$ occuring on its boundary, it follows
that the interior contains a disc $D'$ with circular
$a$-track. Repeating the argument we find that $D$ must have
infinitely many cells which is a contradiction.
\end{proof}

\begin{cor}\label{cor:b-disc} If $D$ is a disc tiled by finitely many
  $[a,b^n]$-discs, then the cyclic word read around $\bdy D$ cannot
  contain only the letter $b$.
\end{cor}

\begin{cor}\label{cor:b-sphere} We cannot tile a sphere with finitely
  many coherently oriented $[a,b^n]$-discs.
\end{cor}

\begin{prop}\label{prop:decomp} Suppose that $D$ is a disc with boundary label
  $[a^N,b^B]$ that is covered by $[a,b^n]$-discs, then it is obtained
  from a collection of $M$ $[a,b^B]$-ribbons $R_1,\ldots R_M$ such that
  the bottom of $R_{i+1}$ is glued to the top of $R_i$, $i=1,\ldots
  M$.
\end{prop}
\begin{proof} We procede by induction on $N$. If $N=1$, then we
  consider the $a-$track $t$ starting at one of the edges of $\bdy D$
  labeled $a$. $t$ must touch the other edge labeled $a$ in $\bdy D$.
  Let $R(t)$ be the subset of $D$ consisting of the $[a,b^n]$-discs
  that $t$ intersects. We note that $R(t)$ can be obtained by making
  some identifications in the top and bottom of some $[a,b^R]$-ribbon,
  but $R(t) \subset D$, which means on one hand that if $R(t)$ is not
  simpy connected then some subset of $\bdy R(t)$ is a circle that
  bounds a disc inside $D$, this disc can only have $b$'s in its label
  contradicting Corollary \ref{cor:b-disc}. It follows that $R(t)$ is
  a ribbon and it contains every $a-$labled edge in $D$, so we must
  have $R(t)=D$.

  Suppose the hypothesis held for all $L \leq N-1$ and suppose that we
  could read $[a^N,b^B]$ along $\bdy D$. We divide $\bdy D$ into four
  arcs $l_a, t_b, r_a, b_b$ that have labels $a^{-N},b^{-B},a^N,b^B$
  respectively, i.e. the left, top, right and bottom sides. Let $e$ be
  the edge with label $a$ that touches the vertex between $l_a$ and
  $t_b$. Let $t$ be the corresponding $a-$track. Let $R(t)$ be as
  above, since $D \subset \E^2$ it is easy to see that $t$ cannot be a
  line from $l_a$ to $l_a$, therefore $t$ must go from $l_a$ to some
  edge $e'$ in $r_a$. 

  Suppose towards a contradiction that $e'$ was not the edge in $r_a$
  that touched the vertex $v$ between $t_b$ and $r_a$. Let $f$ be the
  edge in $r_a$ that touches $v$, and let $u$ be the corresponding
  $a-$track, since $a$-tracks cannot cross we have that $u$ must also
  end in $r_a$ which is a contradiction.

  By the same argument as in the case $N=1$ we have that $R(t)$ must
  be a embedded ribbon. By Corollary \ref{cor:b-disc} we must have
  that $t_b$ is contained in the top of $R(t)$, which means that
  $R(t)$ is an embedded $[a,b^B]$-ribbon and if we remove $R(t)$ from
  $D$, then what remains is a disc $D'$ such that we can read
  $[a^{N-1},b^B]$ along the boundary. So by induction the result
  follows.
\end{proof}

\subsection{A special genus zero quadratic equation}
Equipped with Proposition \ref{prop:decomp} we shall deduce NP
hardness of the following equation:\bel{eqn:quad-bin}
\prod_{j=1}^{k}z_{j}^\mo[a,b^{n_j}]z_j = [a^N,b^B]\ee By the results
in section \ref{sec:olshanskii}, (\ref{eqn:quad-bin}) has a solution if
and only if there is a collection of discs $D_j$ with boundary labels
$[a,b^{n_j}]$ for $j=1\ldots k$ respectively and a disc $D_m$ with
boundary label $[a^N,b^B]$ such that, glued together in a way that
respect labels and orientation of edges, form a union of spheres (this
is forced by the first inequality in (iii), section \ref{sec:olshanskii}).

\begin{thm}\label{thm:np-hard}Deciding if the quadratic equation
  (\ref{eqn:quad-bin}) with coefficients
  \[[a,b^{n_1}],\ldots,[a,b^{n_k}]\tr{~and~} [a^N,b^B]\] has a
  solution is equivalent to deciding if problem
  \ref{prob:mod-bin-packing}; with input $(n_1,\ldots,n_m)$ and
  positive integers $B,N$; has a positive answer. 
\end{thm}
\begin{proof} ``Bin packing $\Rightarrow$ solution.'' Suppose that
  Problem \ref{prob:mod-bin-packing} has a positive answer on the
  specified inputs. For each subset $S_i$ of the given partition of
  $\{1,\ldots,k\}$ we form a $[a,b^B]$-ribbon $R_i$ by glueing
  together the $[a,b^{n_j}]$-discs for $j \in S_i$, this is possible
  by (iv) in section \ref{sec:olshanskii} and equation
  (\ref{eqn:mod-bin}). We then construct one hemisphere by glueing the
  ribbons $R_1,\ldots, R_N$. The other hemisphere is the remaining disc
  with boundary label $[a^N,b^B]^\mo$, the resulting sphere proves the
  solvability of (\ref{eqn:quad-bin}) with the given coefficients.

  ``Solution $\Rightarrow$ bin packing.'' If (\ref{eqn:quad-bin}) has
  a solution then there is a union of spheres tiled with
  $[a,b^{n_i}]$-discs and one $[a^N,b^B]^\mo$-disc, moreover these
  discs are coherently oriented. By condition (v) and Corollary
  \ref{cor:b-sphere} there can only be one sphere: the sphere $S_0$
  containing the unique $[a^N,b^B]^\mo$-disc.  If we remove this
  $[a^N,b^B]^\mo$-disc from $S_0$ what remains will be a disc $D$ with
  boundary label $[a^N,b^B]$ tiled with $[a,b^{n_i}]$-discs. Applying
  Proposition \ref{prop:decomp} divides $D$ into ribbons $R_1,\ldots
  R_N$ and we immmediately see that these ribbons provide a partition
  of $\{n_1,\ldots n_k\}$, showing that Problem
  \ref{prob:mod-bin-packing} has a positive solution on the given
  input.
\end{proof}

\bibliographystyle{alpha} \bibliography{biblio.bib}
\end{document}